\documentclass[11pt, reqno]{amsart}
\usepackage{geometry}                
\usepackage{graphicx}
\usepackage{amssymb}
\usepackage[mathscr]{eucal}
\usepackage[latin1]{inputenc}
\usepackage[brazil]{babel}
\usepackage{hyperref, color}

\newtheorem*{example*}{Exemplo}

\def\R{\ensuremath {\mathbb R}}

\title[O que esperar do calor em uma placa estreita ?]{O que esperar do calor em uma placa estreita ?}

\author[V. M. Nascimento]{Vanderlei M. Nascimento}
\address[Vanderlei Marcos do Nascimento]{Instituto de Geoci\^{e}ncias e Ci\^{e}ncias Exatas, Universidade Estadual Paulista, CEP 13506-900, Rio Claro SP, Brazil}
\email{vandermn@rc.unesp.br}

\author[R. P. Silva]{Ricardo P. Silva$^\star$}
\thanks{$^\star$Apoio FAPESP 2008/53094-4, PROPe-UNESP}
\address[Ricardo Parreira da Silva]{Instituto de Geoci\^{e}ncias e Ci\^{e}ncias Exatas, Universidade Estadual Paulista, CEP 13506-900, Rio Claro SP, Brazil}
\email{rpsilva@rc.unesp.br}

\begin{document}

\begin{abstract}
Ilustramos com um simples exemplo de perturbação de domínios, como a decomposição espectral de um operador linear pode ser aplicada para determinar o comportamento de uma família de equações diferenciais parciais.
\end{abstract}

\maketitle


\section{Introdu\c cão}\label{sec:intro}

Na modelagem matemática de fenômenos naturais, é comum presumir que certas grandezas são desprezíveis comparadas à outras. Por exemplo a profundidade de um oceano comparado à sua extensão. De todo modo, a grandeza desprezível existe tanto quanto as outras, e a matemática leva isso em conta, através da noção de limite: Se $\epsilon$ existe no modelo, o que podemos dizer ao fazermos $\epsilon \to 0$ ? Em cada caso, explicar o significado que se dá para ``quando fazemos $\epsilon \to 0$'', é explicar o que será entendido por ``desprezar $\epsilon $ ''.

Muitos desses aspectos estão incorporados nos chamados problemas em domínios finos, sobre os quais  \cite{Raugel} contempla um panorama substancial.

O que consideramos neste trabalho exemplifica algumas idéias envolvidas no assunto, e, por isso mesmo, apóia-se em resultados importantes já estabelecidos. De posse deles, o tratamento que damos ao exemplo vem com o objetivo de valorizar ferramentas básicas da Álgebra Linear no tratamento de EDP's não lineares. Elas, juntamente com  o fato de podermos tornar preciso o significado de ``fazer $\epsilon \to 0$'', neste caso, permitirão darmos uma resposta quanto ao que esperamos ao fazer isso.

\section{Motivação}\label{sec:motiv}

Imagine que queiramos estudar como o calor se propaga em placas retangulares de espessuras que vão tornando-se mais e mais finas. Fixando o comprimento das placas constante, digamos igual a $1$, como protótipo tomamos o retângulo $\Omega^\epsilon:= (0,1) \times (0,\epsilon)$, $\epsilon >0$, e consideramos $\epsilon \to 0$.

Além disso, para um modelo mais simplificado, supomos que as placas são compostas de material homogêneo e estão termicamente isoladas. Isso significa que o fluxo de calor, sem fontes externas, e com distribuição inicial de temperatura, $u_0^\epsilon$, supostamente  conhecida em cada ponto $(x,y)$ do retângulo $\Omega^\epsilon$, é descrito pelas equações 
\begin{equation}\tag{$\tilde P_\epsilon$}\label{eq:calor-pert}
\left\{
\begin{gathered}
u^\epsilon_t = u^\epsilon_{xx} + u^\epsilon_{yy} , \text{ em }   \R^+ \times \Omega^\epsilon  \\
\frac{\partial u^\epsilon}{ \partial \eta^\epsilon}= 0, \text{ em }  \R^+ \times \partial  \Omega^\epsilon \\
u^\epsilon(0,\cdot,\cdot) = u_0^\epsilon,  \text{ em }   \Omega^\epsilon,
\end{gathered}
\right.
\end{equation}
onde $\eta^\epsilon$ denota o campo vetorial (unitário) normal às laterais $\partial  \Omega^\epsilon$ da placa $\Omega^\epsilon$ e $u^\epsilon=u^\epsilon(t,x,y)$ representa a temperatura no ponto $(x,y)$ de $\Omega^\epsilon$ no instante $t > 0$. Observe que $\dfrac{\partial u^\epsilon}{ \partial \eta^\epsilon}(t,x,y)=~0$ nos pontos de $\partial \Omega^\epsilon$, nos diz que não há quaisquer perda ou ganho de calor através das laterais $\partial  \Omega^\epsilon$, em qualquer instante $t > 0$.

O operador de Laplace
$$
\Delta u (x,y)= u_{xx} (x,y) + u_{yy}(x,y),
$$ 
é um operador linear, e é na análise de seus autovalores e autofunções que seremos capazes de determinar o comportamento limite do fluxo de calor em $\Omega^\epsilon$, quando $\epsilon \to 0$. Note que, quando $\epsilon \to 0$, $\Omega^\epsilon$ se colapsa no segmento $(0,1)$, e assim suspeitamos que se tal comportamento limite existir, deverá ser regido por um modelo unidimensional.

\section{Análise Espectral}\label{espectral}

\subsection{Problema de autovalor}\label{subsec:auto-valor}

Para $\Omega^\epsilon = (0,1) \times (0,\epsilon)$ consideremos o problema de autovalor
\begin{equation}\tag{$\tilde{A}_\epsilon$}\label{eq:av-lap-thin}
\left\{
\begin{gathered}
- \Delta u = \lambda u, \text{ em }   \Omega^\epsilon  \\
\frac{\partial u}{ \partial \eta^\epsilon}= 0, \text{ em }   \partial  \Omega^\epsilon.
\end{gathered}
\right.
\end{equation}

Para cada $\epsilon >0$ fixo, através do método de separa\c cão de variáveis (veja \cite{Djairo}), é fácil obter que (todos) os valores de $\lambda$ para os quais existe uma fun\c cão não nula, $u=u(x,y)$, que satisfaz a equa\c cão \eqref{eq:av-lap-thin},  são 
$$
\lambda^\epsilon_{m,n} = \pi^2(m^2 + \dfrac{n^2}{\epsilon^2}), \quad m,n = 0, 1,2, \cdots
$$ 
 e que, além disso, a respectiva autofunção associada a $\lambda^\epsilon_{m,n}$ é
$$
{u}^\epsilon_{m,n}(x,y) = \cos(m \, \pi\,  x)  \cos(\dfrac{n \, \pi}{\epsilon} \, y).
$$

Munido da norma
$$
\|u\|= \left( \int_{\Omega^\epsilon} |u(x,y)|^2 \, dxdy \right)^{\frac{1}{2}},
$$
 o espaço das (classes de equivalência das) funções cujo quadrado é integrável a Lebesgue (veja por exemplo \cite{Royden}), é completo e usualmente denotado por $L^2(\Omega^\epsilon)$. Também, como esta norma é induzida pelo produto interno
$$
\langle u, v \rangle = \int_{\Omega^\epsilon} u(x,y) v(x,y)  \, dxdy,
$$
$L^2(\Omega^\epsilon)$ é um espaço de Hilbert. Autofunções, $u^\epsilon_{m,n}$ e $u^\epsilon_{m',n'}$, associadas a autovalores distintos, $\lambda^\epsilon_{m,n}$ e $\lambda^\epsilon_{m',n'}$ respectivamente, são ortogonais. De fato, integrando por partes, vê-se, que para quaisquer $(m,n) \neq (m',n') $, 
\begin{align*}
\langle u^\epsilon_{m,n}, u^\epsilon_{m',n'} \rangle & =  \int_{\Omega^\epsilon} \cos(m \, \pi\,  x)  \cos(\dfrac{n \, \pi}{\epsilon} \, y) \, \cos(m' \, \pi\,  x)  \cos(\dfrac{n' \, \pi}{\epsilon} \, y) \, dxdy \\ 
& =  \int_{0}^1 \cos(m \, \pi\,  x) \cos(m' \, \pi\,  x) \, dx \;  \int_{0}^\epsilon  \cos(\dfrac{n \, \pi}{\epsilon} \, y) \cos(\dfrac{n' \, \pi}{\epsilon} \, y) \, dy = 0.
\end{align*}
Além disso,
$$
\| {u}^\epsilon_{m,n} \|=\left[ \int_{\Omega^\epsilon} |{u}^\epsilon_{m,n}|^2  d x d y \right]^{\frac{1}{2}} = \left[  \int_0^1 |\cos(m \, \pi\,  x)|^2 dx \, \int_0^\epsilon |\cos(\dfrac{n \, \pi}{\epsilon} \, y)|^2 dy \right]^{\frac{1}{2}}= \frac{\sqrt{\epsilon}}{a_{m,n}},
$$
onde $a_{0,0}=1$, $a_{m,0}=a_{0,n} = \sqrt{2}$ e $a_{m,n}=2$ para $m,n=1,2,\cdots$. 

Portanto as funções
$$
{w}^\epsilon_{m,n}(x,y) = \dfrac{a_{m,n}}{\sqrt{\epsilon}} \cos(m \, \pi\,  x)  \cos(\dfrac{n \, \pi}{\epsilon} \, y), \quad m,n=0,1,2,\cdots
$$
constituem um sistema ortonormal em $L^2(\Omega^\epsilon)$. Como também é possível mostrar, (veja \cite{Brezis}), que esse sistema é completo temos 
\begin{equation}\label{eq:serie}
u = \sum_{m,n=0}^\infty \langle u, {w}^\epsilon_{m,n}\rangle \, {w}^\epsilon_{m,n}, \quad \forall u \, \in L^2(\Omega^\epsilon).
\end{equation}
%

É fácil ver para $T:L^2(\Omega^\epsilon) \to L^2(\Omega^\epsilon)$ linear e contínuo, que a convergência das somas parciais em \eqref{eq:serie} implica em  $Tu = \displaystyle \sum_{m,n=0}^\infty \langle u, {w}^\epsilon_{m,n}\rangle \, T{w}^\epsilon_{m,n}$. 

Entretanto, mesmo não sendo $-\Delta$ um operador contínuo de $L^2(\Omega^\epsilon)$, também temos\footnote{Ao leitor interessado sugerimos a seção sobre operadores compactos e auto-adjuntos em \cite{Brezis}.}
\begin{equation}\label{eq:serie-Delta}
-\Delta u =\sum_{m,n=0}^\infty\langle u, {w}^\epsilon_{m,n}\rangle \,  \lambda^\epsilon_{m,n} \, {w}^\epsilon_{m,n}, \quad \forall \, u \in {\rm{Dom}}(-\Delta),
\end{equation}
onde ${\rm{Dom}}(-\Delta) = \{u \in L^2(\Omega^\epsilon) : \displaystyle \sum_{m,n=0}^\infty [\langle u, {w}^\epsilon_{m,n} \rangle \,  \lambda^\epsilon_{m,n} ]^2 < \infty \}$.


\subsection{Mudando de variáveis}
  
Ao considerarmos a mudança de variáveis 
$$
\Phi_\epsilon:(x_1,x_2) \mapsto (x, y):=(x_1,\epsilon x_2)
$$ 
que transforma a placa (fixa) $\Omega=(0,1) \times (0,1)$ na placa $\Omega^\epsilon$, se $u=u(x,y)$ é solução do problema \eqref{eq:av-lap-thin}, então, pela regra da cadeia, $v = u \circ \Phi_\epsilon$ é solução do problema
\begin{equation}\tag{$A_\epsilon$}\label{eq:av-lap-thin-rees}
\left\{
\begin{gathered}
- v_{x_1x_1} - \frac{1}{\epsilon^2} v_{x_2 x_2} = \lambda v, \text{ em }    \Omega  \\
v_{x_1} \eta_{x_1} + \frac{1}{\epsilon^2} v_{x_2} \eta_{x_2}= 0, \text{ em }    \partial  \Omega ,
\end{gathered}
\right.
\end{equation}
onde $\eta=(\eta_{x_1},\eta_{x_2})$ denota o campo vetorial (unitário) normal às laterais $\partial  \Omega$ da placa $\Omega$. Similarmente, se $v=v(x_1,x_2)$ é solução de \eqref{eq:av-lap-thin-rees}, então $u = v \circ \Phi_\epsilon^{-1}$ é solução de \eqref{eq:av-lap-thin}. Portanto a mudança $\Phi_\epsilon$ transforma soluções do problema \eqref{eq:av-lap-thin} em soluções do problema \eqref{eq:av-lap-thin-rees}. Da mesma forma, as soluções do problema \eqref{eq:calor-pert} são transformadas em soluções do problema
\begin{equation}\tag{$ P_\epsilon$}\label{eq:calor-mud-var}
\left\{
\begin{gathered}
v^\epsilon_t = v^\epsilon_{x_1x_1} + \frac{1}{\epsilon^2} v^\epsilon_{x_2x_2}  , \text{ em }   \R^+ \times \Omega  \\
v^\epsilon_{x_1} \eta_{x_1} + \frac{1}{\epsilon^2} v^\epsilon_{x_2} \eta_{x_2}= 0, \text{ em }   \partial  \Omega , \\
v^\epsilon(0,\cdot,\cdot) = v_0^\epsilon,  \text{ em }   \Omega,
\end{gathered}
\right.
\end{equation}
considerando-se $v_0^\epsilon = u_0^\epsilon \circ \Phi_\epsilon$.

Cálculos diretos mostram que, para cada $\epsilon >0$, os autovalores do operador 
$$
-\Delta_\epsilon v  = -v_{x_1x_1} - \frac{1}{\epsilon^2} \, v_{x_2x_2}
$$
são os mesmos $\lambda^\epsilon_{m,n} =  \pi^2(m^2 + \dfrac{n^2}{\epsilon^2})$, $m,n=0,1,2,\cdots$, porém as respectivas autofunções (normalizadas) agora independem de $\epsilon$ e são dadas por
$$
v^\epsilon_{m,n}(x_1,x_2) = a_{m,n} \, \cos(m \, \pi\,  x_1)  \cos(n \, \pi \, x_2), \quad m,n=0,1,2,\cdots .
$$

Como antes, as funções $v^\epsilon_{m,n}$, $m,n=0,1,2 \cdots$, constituem um sistema ortonormal completo de $L^2(\Omega)$, e
$$
-\Delta_\epsilon u =  \sum_{m,n=0}^\infty\langle u, {v}^\epsilon_{m,n}\rangle \,  \lambda^\epsilon_{m,n}  \, {v}^\epsilon_{m,n}, \quad \forall \, u \in {\rm{Dom}}(-\Delta_\epsilon).
$$

Ainda por separação de variáveis, somos capazes de expressar a solução $v^\epsilon$ do problema \eqref{eq:calor-mud-var} como 
\begin{equation}\label{eq:sol-calor}
v^\epsilon(t,x_1,x_2) = \sum_{m,n=0}^\infty e^{- t \lambda^\epsilon_{m,n}} \, \langle v_0^\epsilon, {v}^\epsilon_{m,n}\rangle \,     \, {v}^\epsilon_{m,n}(x_1,x_2) .
\end{equation}

Escrevamos agora, para cada $\epsilon >0$, os autovalores $\lambda^\epsilon_{m,n}$ em ordem crescente (levando em conta suas multiplicidades 
$$
\lambda^\epsilon_1 \;  \le \; \lambda^\epsilon_2  \; \le  \; \lambda^\epsilon_3 \; \le  \; \cdots .
$$

Como, para cada natural $k$, existe um valor de parâmetro $\epsilon_k$, que pode ser tomado de modo a que $\pi^2 k^2 \le \dfrac{\pi^2}{\epsilon_k^2}$, então, sempre que $0<\epsilon \le \epsilon_k$, 
$$
\lambda^\epsilon_1 =  0 , \quad \lambda^\epsilon_2 = \lambda^\epsilon_{1,0}, \quad \lambda^\epsilon_3 = \lambda^\epsilon_{2,0}, \quad \cdots,  \quad \lambda^\epsilon_{k+1} = \lambda^\epsilon_{k,0}.
$$

Portanto
$$
\lambda^\epsilon_n  \stackrel{\epsilon \to 0}{\longrightarrow} \pi^2 (n-1)^2 := \lambda^0_{n},
$$  
para todo $n=1,2,\cdots$. 

Da mesma forma, a sequência de autofunções $v^\epsilon_{n}$, associada à sequência $\lambda^\epsilon_n$,  converge trivial\-mente (quando $\epsilon \to 0$),  para a função $v^0_{n}(x_1,x_2): = \sqrt{2} \, \cos(n \, \pi\,  x_1), \, \forall \; n =1,2,\cdots$.

Este comportamento dos autovalores do operador $-\Delta_\epsilon$, no qual o primeiro autovalor ``associado à decomposição na direção $x_2$'' explode quando $\epsilon \to 0$, é fundamental para a caracterização do limite para \eqref{eq:calor-mud-var}, pois ao reconhecermos $\lambda^0_n$ e  $v^0_{n}$, $n=1,2,\cdots$, observamos que são, respectivamente, os autovalores e autofunções (normalizadas) do seguinte problema unidimensional
\begin{equation}\tag{$A_0$}\label{eq:av-lim}
\left\{
\begin{gathered}
- v_{xx}  = \lambda v, \text{ em }    (0,1)  \\
v_x(0) = v_x(1) = 0.
\end{gathered}
\right.
\end{equation}

Similarmente, o espaço $L^2(0,1)$ das funções definidas no intervalo $(0,1)$, com quadrado integrável a Lebesgue, e munido do produto interno
$$
( u, v ) = \int_0^1 u(x) v(x) \, dx
$$
é um espaço de Hilbert, as funções $v^0_{n}$ constituem um sistema ortonormal completo, e
$$
-\frac{d^2}{dx^2} \, v = \sum_{k=1}^\infty  \, ( v, {v}^0_{k} ) \,  \lambda^0_k \, {v}^0_{k}, \quad \forall \, v \in {\rm{Dom}}(-\frac{d^2}{dx^2}),
$$
onde ${\rm{Dom}}(\displaystyle -\frac{d^2}{dx^2}) = \{v \in L^2(0,1) : \displaystyle \sum_{k=1}^\infty \left[ ( v, {v}^0_{k}) \,   \lambda^0_k \right]^2 < \infty \}$.

Também a solução $v^0=v^0(t,x)$ do problema unidimensional
 \begin{equation}\tag{$P_0$}\label{eq:calor-limite}
\left\{
\begin{gathered}
v_t = v_{xx} , \text{ em }   \R^+ \times (0,1) \\
v_x(t,0) = v_x(t,1) =0, \; t > 0\\
v(0,\cdot) = v_0^0,  \text{ em }    (0,1) .
\end{gathered}
\right.
\end{equation}
é dada na forma
\begin{equation}
v^0(t,x) = \sum_{k=1}^\infty e^{- t \lambda^0_{k}} (v_0^0, {v}^0_{k})\,   {v}^0_{k}(x) .
\end{equation}


\section{Redução da dimensão}


Dada uma função $u: (0,1) \to \R$, consideramos $\hat u: \Omega \to \R$ definida por $\hat u(x,y) = u(x)$. Assim  podemos ver o espaço $L^2(0,1)$ como um subespaço de $L^2(\Omega)$, e considerar o operador $-\dfrac{d^2}{dx^2}$ atuando neste subespaço. Na sequência usaremos sem qualquer menção esta identificação.

Com a notação acima, os vetores $\hat v^0_{n}$ constituem um sistema ortonormal completo de  $L^2(0,1)$, e 
$$
\hat v^0(t,x,y) = \sum_{k=1}^\infty e^{- t \lambda^0_{k}}  \langle \hat v_0^0, \hat{v}^0_{k}\rangle \, \hat{v}^0_{k}(x,y) .
$$

Portanto, da convergência espectral obtida na Seção \ref{espectral}, se a sequência de dados iniciais $v_0^\epsilon \in L^2(\Omega)$, converge para $\hat v_0^0$ em $L^2(\Omega)$, então

\begin{equation}
v^\epsilon(t,\cdot,\cdot) = \sum_{k=1}^\infty e^{- t \lambda^\epsilon_{k}} \, \langle v_0^\epsilon, {v}^\epsilon_{k}\rangle \,   {v}^\epsilon_{k} \stackrel{\epsilon \to 0}{\longrightarrow}  \sum_{k=1}^\infty e^{- t \lambda^0_{k}} \, \langle \hat v_0^0, \hat{v}^0_{k}\rangle  \, \hat{v}^0_{k} = \hat v^0(t, \cdot, \cdot),
\end{equation}
uniformemente para $t$ em intervalos compactos de $\R^+$. 
Sobre estes fatos indicamos ao leitor as monografias \cite{Engel-Nagel} e \cite{Gomes}.

Portanto, após essa nossa análise, se perguntados fôssemos :

-- O que esperar do calor em uma placa estreita ?

\noindent de pronto responderíamos:  

-- Esperamos que ele se remeta ao caso unidimensional.


\bigskip \par\noindent {\bf Agradecimentos.} Gostaríamos de expressar nossa gratidão aos professores Antônio Luis Pereira e Luís Augusto F. de Oliveira do IME/USP, cuja discussão informal inspirou a escrita deste texto. Agradecemos também à Profa. Karina Schiabel do DM/UFSCar e ao anônimo revisor pelas sugestões que o aprimoraram.  


\end{document}